\documentclass{article}
\usepackage{amsfonts}
\usepackage{amsmath}
\usepackage
{graphicx}

\usepackage{epsfig}
\begin{document}

\title{Order from Randomness}
\author{Krzysztof Ma\'{s}lanka$^{1}$, Jerzy Cis\l o$^{2}$ \\
\\
$^{1}$Institute for the History of Science, Polish Academy of Sciences\\
kmaslank@uj.edu.pl\\
$^{2}$Institute of Theoretical Physics, University of Wroc\l aw\\
cislo@ift.uni.wroc.pl}
\maketitle

\begin{abstract}
We consider an elementary discrete process which starts from purely random
configuration and leads to well-ordered and stable state. Complete
analytical solution to this problem is presented.
\end{abstract}

\section{Statement of the problem}

Simple models of evolution may sometimes lead to complicated and
unexpected behavior, cf. e.g. \cite{Wolf}.
We consider the following elementary problem. Let us choose a finite set of $%
n$ points $\left\{ p_{i}\right\} $\ distributed arbitrarily (in particular
randomly) on the plane. Let us further label these points in a unique but
otherwise quite arbitrary way:%
\begin{equation}
p_{0},\ p_{1},\ldots ,\ p_{n-1}  \label{points}
\end{equation}%
Of course, given any set one has $n!$ different ways of labelling its
elements. We can further treat these points as two-dimensional sequence of
vectors attached to the origin of the Cartesian coordinates.

Define next the following "cyclic" sequence of vector differences and denote
them $p_{i}(1)$:%
\begin{eqnarray}
p_{0}(1) &=&p_{1}-p_{0}  \label{diffs} \\
p_{1}(1) &=&p_{2}-p_{1}  \notag \\
&&\cdots  \notag \\
p_{n-2}(1) &=&p_{n-1}-p_{n-2}  \notag \\
p_{n-1}(1) &=&p_{0}-p_{n-1}  \notag
\end{eqnarray}%
It represents a sequence of $n$ new vectors which, after appropriate
shifting, may also be attached to the origin of coordinates. We can consider
these points as the next phase of some discrete time evolution and repeat
the process iteratively many times. Note that this process is strictly
deterministic, that is, once chosen the initial points and their labelling,
the rest of the discrete evolution is fully determined regardless of how
many steps has been performed. The only chaotic element in this process is
included in the initial distribution, unless this distribution is
deliberately chosen as regular.

Each transition $t\rightarrow t+1$ in this discrete evolution can be written
using more suggestive matrix notation as:%
\begin{equation}
\left(
\begin{array}{c}
x_{0}(t+1) \\
x_{1}(t+1) \\
~\vdots \\
x_{n-2}(t+1) \\
x_{n-1}(t+1)%
\end{array}%
\right) =\left(
\begin{array}{ccccc}
-1 & +1 & 0 & \cdots & 0 \\
0 & -1 & +1 & \cdots & 0 \\
\vdots & \vdots & \vdots & \vdots & \vdots \\
0 & 0 & \cdots & -1 & +1 \\
+1 & 0 & \cdots & 0 & -1%
\end{array}%
\right) \left(
\begin{array}{c}
x_{0}(t) \\
x_{0}(t) \\
\vdots \\
x_{n-2}(t) \\
x_{n-1}(t)%
\end{array}%
\right)  \label{matrix}
\end{equation}%
for the $x$ coordinate and similarly for the $y$ coordinate, where $%
t=0,1,2,...$ Consequently, starting from the initial configuration, after $t$
evolution steps, we get:%
\begin{equation}
\left(
\begin{array}{c}
x_{0}(t) \\
x_{1}(t) \\
\vdots \\
x_{n-2}(t) \\
x_{n-1}(t)%
\end{array}%
\right) =\left(
\begin{array}{ccccc}
-1 & +1 & 0 & \cdots & 0 \\
0 & -1 & +1 & \cdots & 0 \\
\vdots & \vdots & \vdots & \vdots & \vdots \\
0 & 0 & \cdots & -1 & +1 \\
+1 & 0 & \cdots & 0 & -1%
\end{array}%
\right) ^{{\Huge t}}\left(
\begin{array}{c}
x_{0}(0) \\
x_{0}(0) \\
\vdots \\
x_{n-2}(0) \\
x_{n-1}(0)%
\end{array}%
\right)  \label{matrixn}
\end{equation}%
where the power is understood as matrix power, of course. (Matrix which
appears in (\ref{matrix}) and (\ref{matrixn}) is a special case of so-called
circulant matrix, cf. \cite{Weisstein}, for pedagogical review see \cite%
{Gray}.)

Using different, perhaps less evocative notation, this algorithm can also be
written in coordinate representation as:%
\begin{eqnarray}
x_{k}(t) &=&\left( -1\right) ^{t}\sum\limits_{i=0}^{t}\left( -1\right) ^{i}%
\binom{t}{i}x_{k+i}(0)  \label{forwardn} \\
y_{k}(t) &=&\left( -1\right) ^{t}\sum\limits_{i=0}^{t}\left( -1\right) ^{i}%
\binom{t}{i}y_{k+i}(0)  \notag
\end{eqnarray}%
Equations (\ref{forwardn}) coincide with with the higher order differences
\cite{Arfken}.

In the final section complete analytical solution to this problem will be
given using Discrete Fourier Transform (DFT), see e.g. \cite{Arfken}.

\section{Numerical experiments}

This problem may easily be investigated numerically. The experiments
starting from many different random initial distributions give surprising
and decidedly counter-intuitive results. Naively, one would expect that this
algorithm, when starting from purely random distributions, would just
produce more random distributions. However, it turns out that the final
configuration depends simply on whether $n$ is even or odd. In case of $n$
even, after roughly $n$ steps of discrete evolution, we get exactly two
different separate loops of points. One of these loops contains points with
odd labels whereas the other contains points with even labels. On the other
hand, in case of $n$ odd we get a single loop of points. All configurations
are asymptotically point symmetrical with respect to the origin. (In other
words, the center of mass of the system is in the center of coordinates.)

These results may be qualitatively understood. The crucial thing here is
that the number of points $n$\ is finite and they follow a cyclical pattern.
Since the $t$-th order forward derivatives (\ref{forwardn}) involve $t+1$
terms, when $t$ exceeds\ $n$ some numbers $x_{i}(0)$ appear more than once.
When $t\gg n$ then each $x_{i}(0)$ appears many times in a linear yet
complicated combination.

All figures below were made using \textit{Mathematica} (cf. Fig. 5). We also
made some animated gifs which illustrate this phenomenon and are available
from the authors.
\begin{figure}[h*]
{\epsfig{figure=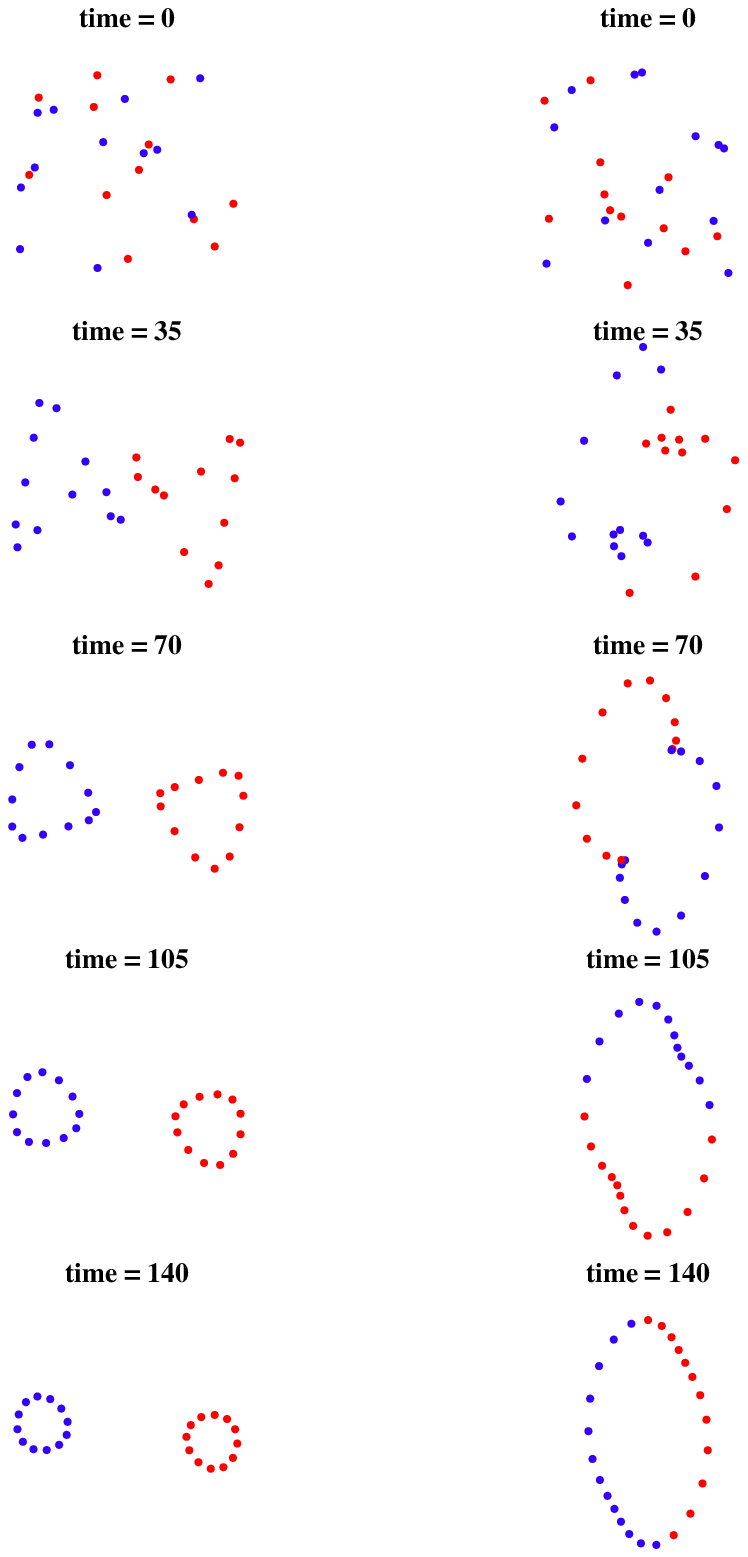, height=15truecm,angle=0}}
\caption{Typical
discrete time evolution for $n=24$ (left) and $n=25$ (right). For even $n$
two different loops emerge whereas for odd $n$ single loop comes out. Since
these loops grow with time, here and in all the remaining figures,
coordinates $x_{i\text{ }}$are normalized, i.e. divided by $\protect\sqrt{%
\sum_{i=0}^{n-1}\left( x_{i}\right) ^{2}}$ (and similarly for~$y$).}
\end{figure}

\begin{figure}[h*]
\centerline{\psfig{file=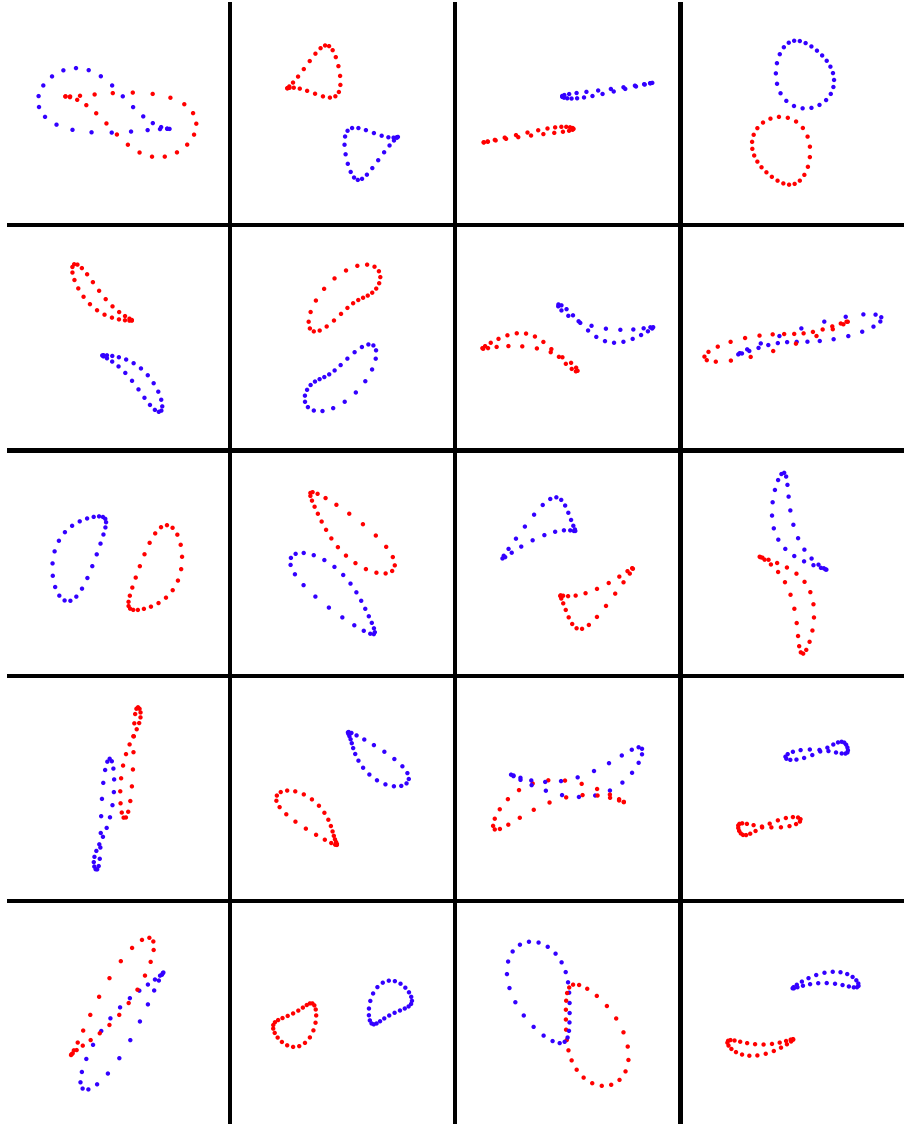,height=15truecm,angle=0}}
\caption{20 cases of initially randomly distributed $n=50$ points
after 300 steps of discrete time evolution. Points $p_{k}$ with odd and even
labels $k$ ($k=0,1,2,...,n-1$) are marked with different colors in order to
show that in the final configuration they belong to different loops and do
not mix. At this stage there are various symmetric shapes but they all tend
to two ellipses, see the next figure.}

\end{figure}

\begin{figure}[h*]
\centerline{\psfig{file=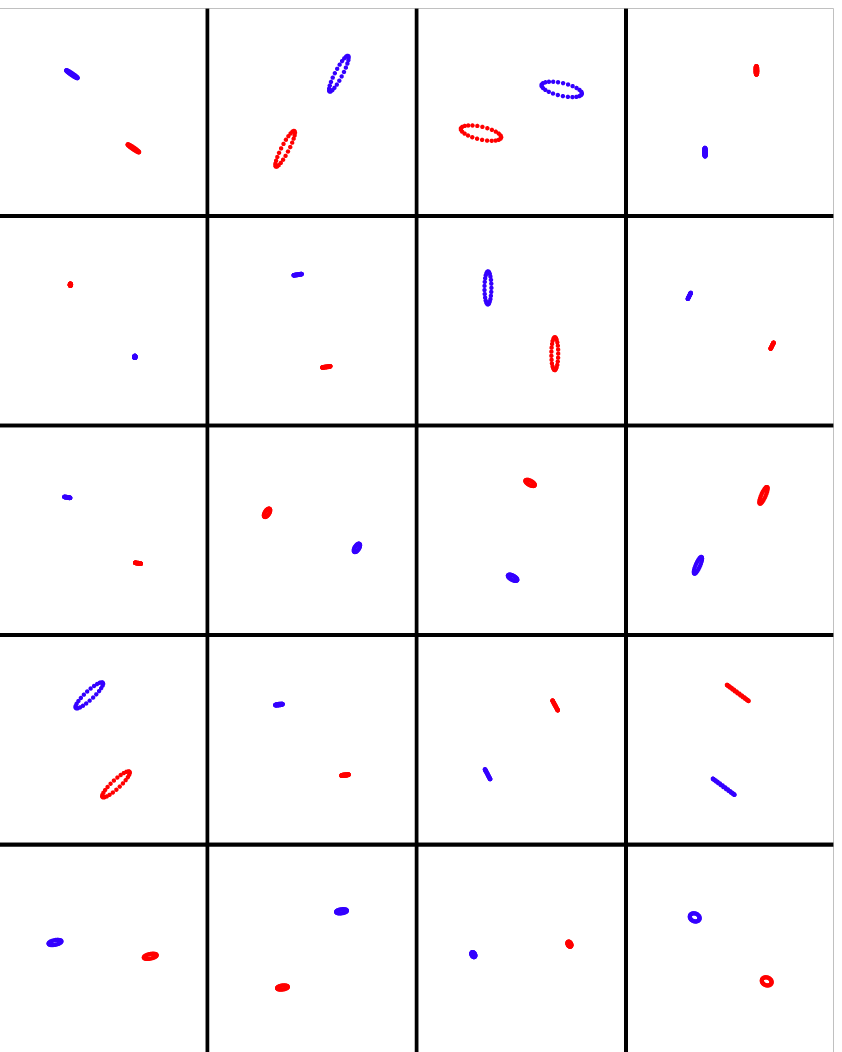,height=15truecm,angle=0}}
\caption{The same as in the
        previous figure but after 1500 steps of discrete time.}
\end{figure}

\begin{figure}[h*]
\centerline{\psfig{file=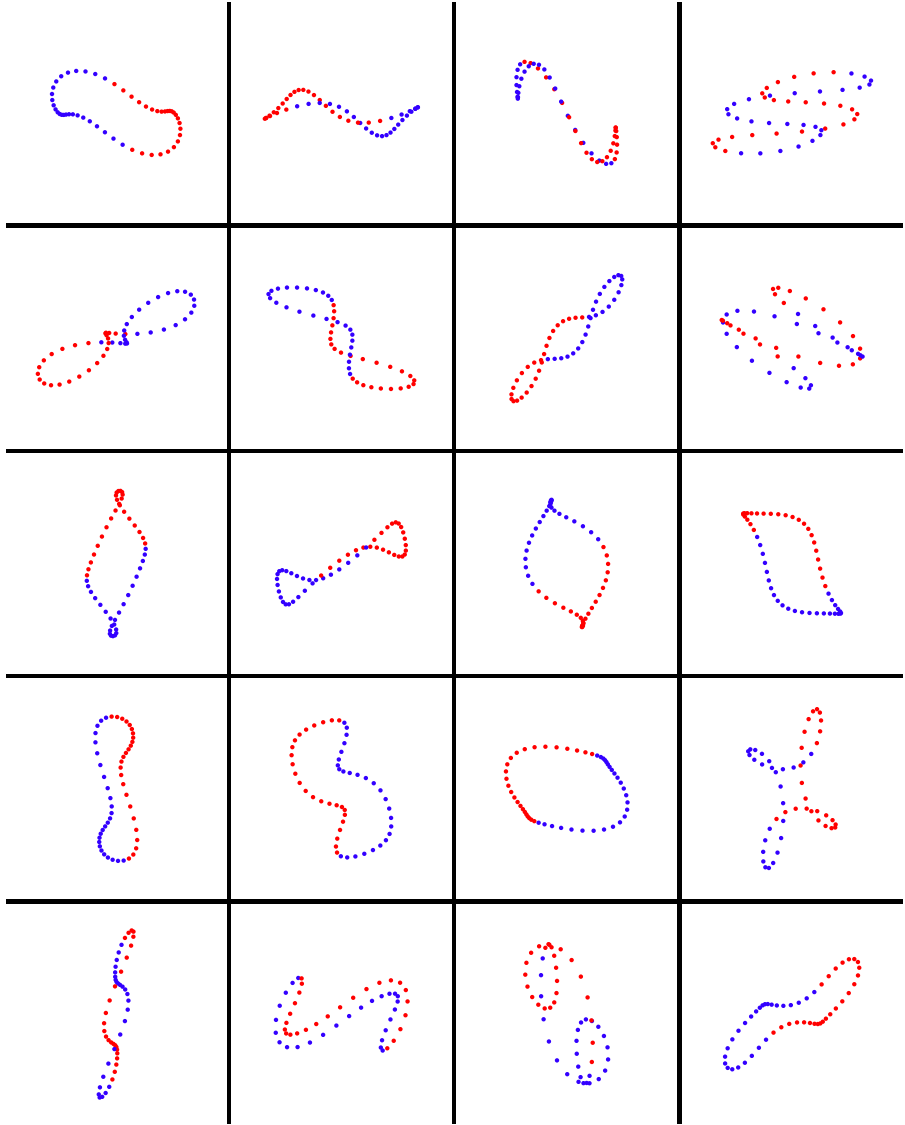,height=15truecm,angle=0}}
\caption{20 cases of initially
        randomly distributed $n=$51 points after 300 steps of discrete time
        evolution. Single loop is visible with various degree of entanglement. Note
        that points with even and odd values of the labelling parameter also do not
        mix. There are many different shapes but they all tend to some ellipse, cf.
        eq. (\protect\ref{ellipse}), see also the next figure.}
\end{figure}

\begin{figure}[h*]
\centerline{\psfig{file=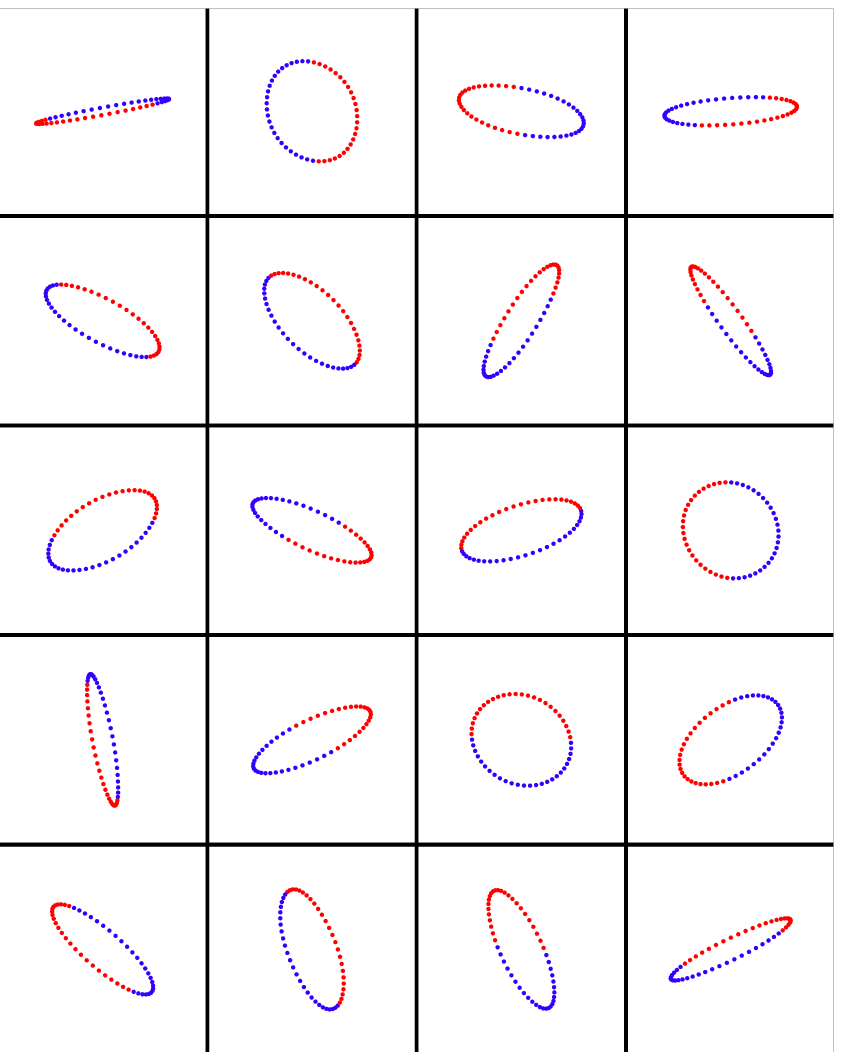,height=15truecm,angle=0}}
\caption{The same as in the
        previous figure but after 1500 steps of discrete time.}
\end{figure}

\begin{figure}[h*]
\centerline{\psfig{file=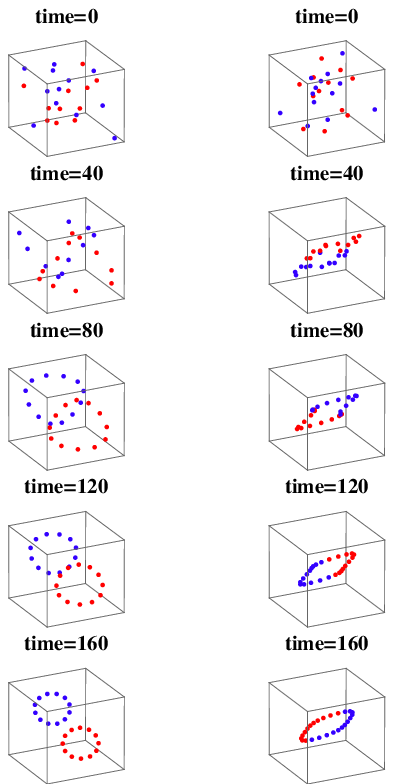,height=15truecm,angle=0}}
\caption{3D case reveals the same qualitative behavior. $n=24$
(left), $n=25$ (right).}
\end{figure}

\section{Analytical solution}

The statement of the problem was given by the first author. He also
performed the numerical experiments described above. The complete analytical
solution which we shall present below was given by the second author.

We may consider the coordinates separately. The calculations below concern
the $x$ coordinate but may be applied to $y$ (as well as to other possible
coordinates in higher dimensions).

Let us denote initial values of coordinates by
\begin{equation*}
x_{0}(0),x_{1}(0),\dots ,x_{n-1}(0).
\end{equation*}%
Recall that the discrete time evolution is described by the equations:%
\begin{equation}
\begin{array}{l}
x_{j}(t+1)=x_{j+1}(t)-x_{j}(t),\qquad t=0,1,\dots ,n-2 \\
x_{n-1}(t+1)=x_{0}(t)-x_{n-1}(t).%
\end{array}
\label{r}
\end{equation}%
The evolution is linear, therefore it may be diagonalized using discrete
Fourier transform
\begin{equation}
\hat{x}_{k}(t)=\frac{1}{\sqrt{n}}\sum_{j=0}^{n-1}x_{j}(t)\omega
^{-jk},\qquad k=0,1,\dots ,n-1,  \label{ft}
\end{equation}%
where $\omega =\exp (2\pi i/n)$. Fourier transform of the eq. (\ref{r})
reads
\begin{equation}
\hat{x}_{k}(t+1)=(\omega ^{k}-1)\hat{x}_{k}(t).  \label{tfr}
\end{equation}%
The general formula is
\begin{equation}
\hat{x}_{k}(t)=(\omega ^{k}-1)^{t}\hat{x}_{k}(0)  \label{gf}
\end{equation}%
The inverse Fourier transform of eq. (\ref{gf}) is
\begin{eqnarray}
x_{l}(t) &=&\frac{1}{\sqrt{n}}\sum_{k=0}^{n-1}\hat{x}_{k}(t)\omega ^{lk}=%
\frac{1}{n}\sum_{k=0}^{n-1}(\omega ^{k}-1)^{t}\omega
^{lk}\sum_{j=0}^{n-1}x_{j}(0)\omega ^{-jk} \\
&=&\frac{1}{n}\sum_{k=0}^{n-1}(\omega ^{k}-1)^{t}\omega
^{lk}\sum_{j=0}^{n-1}x_{j}(0)\omega ^{-jk}\qquad l=0,1,\dots ,n-1.  \notag
\end{eqnarray}%
From this point we will write equations for two coordinates.

In the case of odd $n$, absolute value of $|\omega ^{k}-1|$ reaches
maximum of $r=2\cos (\pi /2n)$ for $k=(n\pm 1)/2$,
 $\ \omega
^{k}=-\exp (\pm \pi i/n)$.
 Hence for large $t$ we have the following
approximation
\begin{equation}
x_{l}(t)\approx r^{t}(-1)^{l+t}\left( \cos \frac{\pi (l+t/2)}{n}\;A+\sin
\frac{\pi (l+t/2)}{n}\;B\right)
\end{equation}%
\begin{equation*}
y_{l}(t)\approx r^{t}(-1)^{l+t}\left( \cos \frac{\pi (l+t/2)}{n}\;C+\sin
\frac{\pi (l+t/2)}{n}\;D\right)
\end{equation*}%
where $A,B,C,D$ are numerical coefficients:
\begin{eqnarray}
A &=&\frac{2}{n}\sum_{j=0}^{n-1}(-1)^{j}\cos \frac{\pi j}{n}\;x_{j}(0)
\label{ABCD} \\
B &=&\frac{2}{n}\sum_{j=0}^{n-1}(-1)^{j}\sin \frac{\pi j}{n}\;x_{j}(0)
\notag \\
C &=&\frac{2}{n}\sum_{j=0}^{n-1}(-1)^{j}\cos \frac{\pi j}{n}\;y_{j}(0)
\notag \\
D &=&\frac{2}{n}\sum_{j=0}^{n-1}(-1)^{j}\sin \frac{\pi j}{n}\;y_{j}(0)
\notag
\end{eqnarray}%
We note that the following relations hold%
\begin{eqnarray}
x_{l}(t+2) &=&-r^{2}x_{l+1}(t) \\
y_{l}(t+2) &=&-r^{2}y_{l+1}(t).  \notag
\end{eqnarray}%
We have in this case%
\begin{equation}
\left(
\begin{array}{c}
x_{l}\left( t\right) \\
y_{l}\left( t\right)%
\end{array}%
\right) \approx r^{t}(-1)^{l+t}\left(
\begin{array}{cc}
A & B \\
C & D%
\end{array}%
\right) \left(
\begin{array}{c}
\cos \phi \\
\sin \phi%
\end{array}%
\right)
\end{equation}%
where%
\begin{equation*}
\phi =\frac{\pi (l+t/2)}{n}
\end{equation*}%
Hence for large $t$ and $AD-BC\neq 0$, the points $p_{0}(t),p_{1}(t),\dots
,p_{n-1}(t)$ lie on an ellipse
\begin{equation}
(C^{2}+D^{2})x^{2}-2(AC+BD)xy+(A^{2}+B^{2})y^{2}=(AD-BC)^{2}r^{2t}.
\label{ellipse}
\end{equation}

In the case of even $n$, absolute value of $|\omega ^{k}-1|$ reaches maximum
value equal $2$ for $k=n/2$, $\ \omega ^{k}=-1$ / Hence for large $t$ we
have an approximate relation
\begin{eqnarray}
x_{l}(t) &\approx &2^{t}(-1)^{l+t}A \\
y_{l}(t) &\approx &2^{t}(-1)^{l+t}C,  \notag
\end{eqnarray}%
where
\begin{eqnarray}
A &=&\frac{1}{n}\sum_{j=0}^{n-1}(-1)^{j}x_{j}(0)  \label{AC} \\
C &=&\frac{1}{n}\sum_{j=0}^{n-1}(-1)^{j}y_{j}(0).  \notag
\end{eqnarray}

More detailed investigations, when higher eingenvalue is taken into account,
give the following asymptotics:%
\begin{equation*}
\binom{x_{l}(t)}{y_{l}(t)}\approx (-1)^{l+t}\left[ \binom{A}{C}%
+r_{1}^{t}\left(
\begin{array}{cc}
A_{1} & B_{1} \\
C_{1} & D_{1}%
\end{array}%
\right) \binom{\cos \phi }{\sin \phi }\right] ,\qquad \phi =\frac{2\pi
\left( l+\frac{t}{2}\right) }{n}
\end{equation*}%
where $r_{1}=2\cos (\pi /n)$\ and $A_{1},B_{1},C_{1},D_{1}$ are defined
similarly to (\ref{ABCD}) with $2\pi $ instead of $\pi $.

\section{Conclusion}

The process presented in this paper may easily be generalized to more
spatial dimensions. The method described in the previous section applies
also in this case. Just as in the 2D case the final configuration depends
solely on whether $n$ is odd or even (cf. Fig. 4).

Let us finally note that this example of discrete time evolution
leading from random to ordered distribution complies with several informal
requirements formulated by some distinguished mathematicians when they spoke
about aesthetics or even beauty of their discipline. Namely, given result
should be easily stated, yet counter-intuitive, deeply non-obvious,
surprising. It should have rigorous and elegant proof. Finally -- as G. H.
Hardy used to emphasize - it should be perfectly useless \cite{Hardy}%
.\bigskip

\textbf{Acknowledgement}

The authors would like to thank Prof. Marek Wolf, Institute of the
Theoretical Physics, University of Wroc\l aw, for stimulating discussion and
for performing several independent numerical experiments using \textit{Delphi%
} which confirmed the numerical results presented in this paper.

\medskip

\subsubsection*{Mathematica code used for generating animated gifs illustrating
the effect described in this paper}
\begin{verbatim}
(* Mathematica code *)

(* M random labelled points on the plane *)

M = 25;
ro := Random[Real, {-1,1}, 10]
Table[{ro, ro, ro}, {m, 1, M}];
points = %;

(* differences between consecutive *)
(* vectors generate M new vectors  *)

F[P_] := Table[
  If[k < Length[P], P[[K+1]]-P[[k]],P[[1]]-P[[Length[P]]]],
  {k, 1, Length[P]}
 ]

(* iteration and normalization *)

            Nest[F, points, i]
G[i_] := --------------------------
          Norm[Nest[F, points, i]]

(* color of points: odd - green, even - red *)

b[i_] := If[ OddQ[i], 0.4, 0.0]
\end{verbatim}

\subsubsection*{Mathemaatica code (continued)}
\begin{verbatim}
(*drawing points *)

H[k_] := Module[
  {t0},
  Show[
   Table[
    Graphics3D[{
      Hue[b[k+i]],
      PointSize[r],
      Point[
       {t0[[i]][[1]], t0[[i]][[2]], t0[[i]][[3]]}]
     }],
    {i, 1, Length[t0]}
   ],
   Background -> RGBColor[0.2, 0.2, 0.3],
   AspectRatio -> 1,
   Axces -> False,
   Boxed -> True,
   PlotRange -> {{-R, +R}, {-R, +R}, {-R, +R}},
   ViewPoint -> {100, -50, 50},
   PlotLabel -> StyleForm["times" <> ToString[k], FontSize -> 12]
 ]

 (* size of dots *)
 r = 0.008;
 (* size of picture *)
 R = 0.13;

 t0 = SessionTime[];
 Table[H[V[k]], {k, 1, 510, 1} ]
 SessionTime[] - t0

 Export["animation.gif", %%, "GIF", ImageResolution -> 100]
\end{verbatim}

\end{document}